\documentclass{amsart}

\usepackage{amsmath,amssymb,amsthm,latexsym}

\newtheorem{theorem}{Theorem}
\newcommand{\bt}{\begin{theorem}}
\newcommand{\et}{\end{theorem}}
\newtheorem{lemma}{Lemma}
\newcommand{\bl}{\begin{lemma}}
\newcommand{\el}{\end{lemma}}
\newtheorem{corollary}{Corollary}
\newcommand{\bc}{\begin{corollary}}
\newcommand{\ec}{\end{corollary}}
\newcommand{\bconj}{\begin{conjecture}}
\newcommand{\econj}{\end{conjecture}}
\newtheorem{problem}{Problem}
\newcommand{\bprob}{\begin{problem}}
\newcommand{\eprob}{\end{problem}}
\newcommand{\beq}{\begin{equation}}
\newcommand{\eeq}{\end{equation}}
\newcommand{\benum}{\begin{enumerate}}
\newcommand{\eenum}{\end{enumerate}}
\newcommand{\N}{\ensuremath{ \mathbf N }}
\newcommand{\Z}{\ensuremath{\mathbf Z}}
\newcommand{\Q}{\ensuremath{\mathbf Q}}
\newcommand{\R}{\ensuremath{\mathbf R}}

\newcommand{\mcr}{\ensuremath{ \mathcal R}}

\newcommand{\mba}{\ensuremath{ \mathbf a}}
\newcommand{\mbb}{\ensuremath{ \mathbf b}}

\newcommand{\mbe}{\ensuremath{ \mathbf e}}

\newcommand{\mbo}{\ensuremath{ \mathbf 0}}
\newcommand{\mbu}{\ensuremath{ \mathbf u}}
\newcommand{\mbv}{\ensuremath{ \mathbf v}}

\newcommand{\mbx}{\ensuremath{ \mathbf x}}
\newcommand{\mby}{\ensuremath{ \mathbf y}}

\DeclareMathOperator{\diam}{\text{diam}}

\newcommand{\bmat}{\left(\begin{matrix}}
\newcommand{\emat}{\end{matrix}\right)}
\newcommand{\bsmallmat}{\left(\begin{smallmatrix}}
\newcommand{\esmallmat}{\end{smallmatrix}\right)}

\DeclareMathOperator{\qqand}{\qquad\text{and}\qquad}

\title[Compression  and complexity for sumset sizes]
{Compression and complexity for sumset sizes in additive number theory}
\author{Melvyn B.  Nathanson}
\address{Department of Mathematics\\Lehman College (CUNY)\\Bronx, NY 10468}
\email{melvyn.nathanson@lehman.cuny.edu}

\date{\today}

\thanks{Supported in part by  PSC-CUNY Research Award Program grant 66197-00 54.}

\begin{document}

\begin{abstract}
The study of sums of finite sets of integers has mostly concentrated on sets 
with small sumsets (Freiman's theorem and related work) and on sets with 
large sumsets (Sidon sets and $B_h$-sets).  This paper considers the sets  
${\mathcal R}_{\mathbf Z}(h,k)$ and ${\mathcal R}_{{\mathbf Z}^n}(h,k)$ 
of \emph{all} sizes of $h$-fold sums of sets of $k$ integers or of $k$ lattice points, 
and the geometric and computational complexity of the sets 
${\mathcal R}_{\mathbf Z}(h,k)$ and ${\mathcal R}_{{\mathbf Z}^n}(h,k)$. 
For sumsets $hA$ with large diameter, there is a compression algorithm 
to construct sets $A'$ with $|hA'| = |hA|$ and small diameter. 
\end{abstract}

\maketitle

\section{Geometric and computational complexity}
Let $X$ be an additive abelian semigroup and let $A$ be a nonempty subset 
of $X$.   The $h$-fold sum of $A$ is the set of all sums of $h$ 
not necessarily distinct elements of $A$. 
We study the sizes of $h$-fold sumsets 
of $k$-element subsets of $X$:  
\[
\mcr_X(h,k) = \{ |hA| : A \subseteq X \text{ and } |A| = k\}. 
\]
A central  problem in additive number theory is to compute 
and understand the sumset size set $\mcr_X(h,k)$.   This is  unsolved  
even in the classical cases when $X$ is the semigroup $\N_0$ 
of  nonnegative integers, 
the groups \Z, \Q, and \R\ of integers,  rational numbers, and real numbers, 
and the groups $\Z^n$, $\Q^n$, and $\R^n$ of lattice points with integral, 
 rational, and real coordinates.   
 If sets $A$ an $A'$ are affinely equivalent, that is, if
\[
A' = \lambda \ast A + \mu =  \{\lambda a + \mu: a \in A \}
\]
for real numbers $\lambda \neq 0$ and $\mu$, then $|hA'| = |hA|$. 
It follows that, for all positive integers $h$ and $k$, 
\beq           \label{compression:1-dim}
\mcr_{\N_0 }(h,k) = \mcr_{\Z }(h,k) = \mcr_{\Q }(h,k) 
\eeq
and 
\beq           \label{compression:n-dim}
 \mcr_{\N_0^n }(h,k) =  \mcr_{\Z^n }(h,k) = \mcr_{\Q^n }(h,k). 
\eeq
We have 
\beq           \label{compression:Z-min}
\min \mcr_{\Z}(h,k) = \min \mcr_{\R }(h,k) = hk-h+1
\eeq
and 
\beq           \label{compression:Z-max}
\max \mcr_{\Z}(h,k) = \max \mcr_{\R }(h,k) = \binom{h+k-1}{k-1}. 
\eeq 
In Theorem~\ref{compression:theorem:LatticePoints}, we prove that 
\[
\mcr_{\Z}(h,k) = \mcr_{\Z^n }(h,k) 
\]
for all positive integers $n$.  
Of course,  
 \[
\mcr_{\Z }(h,k) \subseteq  \mcr_{\R }(h,k)
\]
but it is not known if $\mcr_{\Z}(h,k) = \mcr_{\R }(h,k)$ 
for some or all positive integers $h$ and $k$.  

This paper considers the geometric and computational complexity of the sets 
$\mcr_{\Z}(h,k)$ and $\mcr_{\Z^n }(h,k)$.  This problem was introduced 
in~\cite{nath25aa,nath25bb,nath25cc} and basic results can be found there.
For all positive integers $h$ and $k$, we have   
\[
\mcr_{\Z}(h,1) = \{1\} ,\qquad \mcr_{\Z}(h,2) = \{h+1\}, \qquad  \mcr_{\Z}(1,k) = \{k\} 
\]
and, from~\cite{nath25bb}, 
\[
\mcr_{\Z}(2,k) = \left[ 2k-1,\binom{k+1}{2} \right]
\]
where, for all $u,v \in \R$,there is the \emph{integer interval} 
\[
[u,v] = \{n\in \Z: u \leq n \leq v\}.
\]

\bprob
The fundamental problem is to compute $\mcr_{\Z}(h,k)$ for all 
$h \geq 3$ and $k \geq 3$. 
\eprob

The diameter of a nonempty subset $A$ of a metric space  is 
$\diam(A) = \sup\{d(x,y):x,y \in A\}$. 
Let $A = \{a_1, a_2, \ldots, a_{k-1}, a_k\}$ be a set of $k$ nonnegative integers 
with $0= a_1 < a_2 < \cdots < a_{k-1} < a_k$ and $|hA| = t \in \mcr_{\Z}(h,k)$.  
If $a_k > ha_{k-1} + 1$, then the ``compressed'' set 
$A' = \{a_1, a_2, \ldots, a_{k-1}, a'_k\}$ with $a'_k = ha_{k-1}+1$ 
also satisfies $|A'| = k$ and $|hA'| = t$ 
(Theorem~\ref{compression:theorem:Delta-ak}), and 
\[
\diam(A') = ha_{k-1}+1 < a_k =  \diam(A).
\]
For example, if $A = \{0,1,3,20,200\}$, 
then $|A| = 5$ and $|4A| = 65$. 
If $A' = \{0,1,3,20,81\}$, 
then $|A'| = 5$ and $|4A'| = 65$. 
Indeed, for all $b \geq 81$, if $A' = \{0,1,3,20,b\}$, 
then $|A'| = 5$ and $|4A'| = 65$.  However,    if $A' = \{0,1,3,20,80\}$, 
then $|A'| = 5$ and $|4A'| = 64$.

In this example, we reduced the diameter of the set $A$ 
by compressing its largest element.  Consider reducing the diameter of $A$ by 
compressing the gap between ``interior'' elements of $A$.  
For example, the set $A' = \{0,1,3,20,81\}$ satisfies  $|A'| = 5$ and $|4A'| = 65$. 
Shrinking the gap between the elements 3 and 20, we see that the set 
$A'' = \{0,1,3,20-c,81-c\}$ satisfies  $|A''| = 5$ and $|4A''| = 65$ for all 
$c \in [0,7]$, but, with $c=8$ and $A'' = \{0,1,3,12,73\}$, 
we obtain $|4A''| = 64$.  

\bprob
The geometric complexity problem is to describe, for positive integers $h$ 
and $k$,  the shapes of ``compressed'' 
sets of $k$ integers or lattice points with given $h$-fold sumset size $t$.  
\eprob

There are also computational complexity problems for sumsets.   
By~\eqref{compression:Z-min} and~\eqref{compression:Z-max}, 
the set $\mcr_{\Z}(h,k)$ is finite. 
For all $t \in \mcr_{\Z}(h,k)$, there is an integer $N_t$ and a set $A \subseteq [0,N_t-1]$ 
with $|A| = k$ and $|hA| = t$.  
Let $N = \max\{N_t: t \in \mcr_{\Z}(h,k) \}$. 
Then 
\[  
\mcr_{\Z}(h,k) = \{ |hA|: A \subseteq [0,N-1] \text{ and } |A| = k \}.
\]
Let $N(h,k)$ be the least integer $N$ that satisfies this condition. 
For all $h \geq 1$ and $k \geq 1$ we have 
\[
N(h,1) = 1, \qquad N(1,k) = k, \qqand N(h,2) = 2. 
\]
Nathanson~\cite{nath25bb} proved that 
\[
N(2,k) \leq 2^k \qqand N(3,3)=5.  
\]
Theorem~\ref{compression:theorem:N(h,k)} gives the upper bound 
\[
N(h,k) < 4(8h)^{k-1}   
\]
for all $h \geq 3$ and $k \geq 3$. 

\bprob
The computational complexity problem is to determine $N(h,k)$ 
for all $h \geq 3$ and $k \geq 3$. 
\eprob

For related work on sumset sizes, see~\cite{fox-krav-zhan25, hegy96,krav25,nath25aa,peri-roto25}.

\section{Compression theorem for integers}

By abuse of notation, we write 
\[
A = \{a_1 < a_2 < \cdots < a_k\}
\]
if $A = \{a_1,a_2, \ldots, a_k\}$ is a set of $k$ integers and $a_1 < a_2 < \cdots < a_k$.

\bt                                           \label{compression:theorem:Delta-1}
 Let $k \geq 3$ and let $A = \{a_1 <  a_2 < \cdots < a_k\}$ be a set 
of $k$ integers. 
Let $h$ be a positive integer.  
If there exist an integer $j \in [1,k-1]$ and a positive integer $\delta_j$  such that 
\beq                                     \label{compression:Delta-ineq-2}
a_{j+1} - a_j = 1 + \delta_j +  (h-1)\max (a_j - a_1, a_k - a_{j+1})  
\eeq
then the integer $j$ is unique and the strictly increasing sequence 
of $k$ integers defined by   
\beq                                     \label{compression:a-ell}
a'_{\ell} = \begin{cases}
a_{\ell} & \text{if $\ell \in [1,j]$} \\
a_{\ell} - \delta_j  & \text{if $\ell \in [j+1, k]$  } 
\end{cases}
\eeq
 and the set $A' = \{a'_1 <  a'_2 < \cdots < a'_k\}$ 
satisfy  
\beq           \label{compression:hA}
|hA'| = |hA|.
\eeq 
Moreover, 
\beq                                     \label{compression:Delta-ineq}
a'_{j+1} - a'_j =  1+ (h-1)\max (a'_j - a'_1, a'_k - a'_{j+1})   
\eeq 
and 
\beq           \label{compression:diam} 
\diam(A') =  \diam(A) - \delta_j \leq \diam(A) - 1. 
\eeq
\et

\begin{proof}
Suppose there exist distinct integers $j_1 \in [1,k-1]$ and $j_2\in [1,k-1]$ 
that satisfy~\eqref{compression:Delta-ineq-2}.  
If $j_1 < j_2$, then 
\begin{align*}  
a_{j_1+1} - a_{j_1} 
& >  (h-1)\max (a_{j_1} - a_1, a_k - a_{j_1+1}) \\
& \geq  (h-1)(a_k - a_{j_1+1}) \\
&  \geq  a_{j_2+1} - a_{j_2}  \\
& >   (h-1)\max (a_{j_2} - a_1, a_k - a_{j_2+1}) \\
& \geq   (h-1)(a_{j_2} - a_1) \\
&  \geq  a_{j_1+1} - a_{j_1}  
\end{align*}
which is absurd.  Therefore, at most one integer $j \in [1,k-1]$ satisfies~\eqref{compression:Delta-ineq-2}.

Condition~\eqref{compression:Delta-ineq-2} implies 
\begin{align}         \label{compression:bad-j-2}
 (h   - 1  -i) &  (a_j - a_1) + i(a_k - a_{j+1})       \nonumber  \\ 
 & \leq   (h-1)\max (a_j - a_1, a_k - a_{j+1})  \\
 &   < a_{j+1} - a_j - \delta_j     \nonumber
\end{align} 
for all $i \in [0,h-1]$.  

Let 
\[
A = B \cup C
\]
where 
\[
B =  \{a_1 < \cdots < a_j\}
\]
and 
\[
C =  \{a_{j+1}  < \cdots < a_k\}. 
\]
The sumset $hA$ decomposes as follows:     
\[
hA = h(B \cup C ) = \bigcup_{i=0}^h \left(  (h-i)B + iC \right). 
\]
We have 
\begin{align*} 
(h-i)B + iC &  \subseteq [ (h-i)a_1, (h-i)a_j]+ [ia_{j+1}, ia_k] \\ 
& =  [ (h-i)a_1 + ia_{j+1},(h-i)a_j+ ia_k].   
\end{align*} 
For all $i \in [0,h-1]$,
\[
\max \left( (h-i)B + iC \right) =  (h-i)a_j+ ia_k 
\]
and 
\[
\min \left( (h - (i + 1))B + (i+1)C \right) = (h- 1-i )a_1 + (i+1)a_{j+1}.
\]
The $h+1$ sets $(h-i)B + iC$ are pairwise disjoint if, for all $i \in [0,h-1]$,  
\[ 
 (h-i)a_j+ ia_k  <   (h- 1-i )a_1 + (i+1)a_{j+1} 
\] 
or, equivalently, if 
\[
(h-1-i ) (a_j - a_1) + i(a_k - a_{j+1}) < a_{j+1} -a_j .  
\]
This is inequality~\eqref{compression:bad-j-2} and so 
the $h+1$ sets $(h-i)B + iC$ are  pairwise disjoint and   
\beq          \label{compression:delta}
|hA| =  \sum_{i=0}^h \left|  (h-i)B + iC \right|.
\eeq

We have 
\begin{align*} 
a'_j - a'_1 & = a_j - a_1 \\ 
 a'_k - a'_{j+1} & = (a_k - \delta_j) - (a_{j+1}  - \delta_j) = a_k  - a_{j+1} 
 \end{align*} 
 and, for all $\ell \in [1,k-1] \setminus \{j\}$,  
\begin{align*}
a'_{{\ell} +1} - a'_{\ell} & = a_{{\ell} +1}-a_{\ell} \geq 1. 
\end{align*} 
Condition~\eqref{compression:Delta-ineq-2} gives 
\begin{align*}
 a'_{j+1} - a'_j & = (a_{j+1} - \delta_j) - a_j  \\ 
& = 1 + (h-1)\max (a_j - a_1, a_k - a_{j+1}) \\
& = 1 + (h-1)\max (a'_j - a'_1, a'_k - a'_{j+1}).
\end{align*}  
This proves~\eqref{compression:Delta-ineq}. 
Moreover, 
\[
a'_1 < \cdots < a'_{j-1} < a'_j < a'_{j+1}  < \cdots < a'_k 
\]
and the set 
\[
A'  = \{a'_1 < a'_2 < \cdots < a'_k\}  
\] 
satisfies 
\[
\diam(A') = a'_k-a'_1 = a_k - a_1 - \delta_j = \diam(A) - \delta_j \leq \diam(A) -1.
\]
This proves~\eqref{compression:diam}. 

We shall prove that $|hA'| = |hA|$.  Let 
\[  
A'  = B' \cup C'  
\]  
where 
\[
B' =  \{a'_1< \cdots <  a'_j\} =  \{a_1< \cdots <  a_j\} = B  
\]
and 
\[
C' = \{a'_{j+1} < \cdots <  a'_k\} =   \{a_{j+1} - \delta_j  < \cdots <  a_k - \delta_j \} =  C - \delta_j.  
\]
Then  
\[
hA' = h(B' \cup C' ) = \bigcup_{i=0}^h \left(  (h-i)B' + iC' \right). 
\]
We have 
\[
 (h-i)B' + iC' =  (h-i)B + i(C - \delta_j) =  (h-i)B + iC - i\delta_j
\]
and so 
\[
\left|  (h-i)B' + iC' \right| = \left|  (h-i)B + iC \right|.
\]
If the sets $ (h-i)B' + iC'$ are pairwise disjoint for all $i \in [0,h]$, then 
\[ 
|hA'| = \sum_{i=0}^h \left|  (h-i)B' + iC' \right|  = \sum_{i=0}^h \left|  (h-i)B + iC \right| 
= |hA|. 
\]
Observe that 
\begin{align*} 
(h-i)B' + iC' &  \subseteq [(h-i)a_1 ,(h-i)a_j]+ [ia_{j+1} - i \delta_j , ia_k - i \delta_j ] \\ 
& =  [ (h-i)a_1 + ia_{j+1}  - i \delta_j ,(h-i)a_j+ ia_k  - i \delta_j ]. 
\end{align*} 
For all $i \in [0,h-1]$, the sets $(h-i)B' + iC' $ and $(h-i-1)B' + (i+1)C' $ 
are pairwise disjoint if 
\[
(h-i)a_j+ ia_k  - i \delta_j < (h-i-1)a_1 + (i+1)a_{j+1}  - (i+1)\delta_j 
\]
or, equivalently, if 
\[
(h-i-1)(a_j -a_1)+ i(a_k-a_{j+1}) < a_{j+1}-a_j - \delta_j.  
\]
This is exactly inequality~\eqref{compression:bad-j-2} and so $|hA'| = |hA|$.  
This  completes the proof.  
\end{proof}

\bt                                               \label{compression:theorem:Delta-ak}
Let $A = \{a_1 < \cdots < a_{k-1} < a_k\}$ be a set of $k$ integers  with $a_1=0$.  
For all positive integers $h$, if $a_k > 1 + ha_{k-1}$, 
then the set $A' = \{a_1 < \cdots < a_{k-1} < a'_k \}$
with $a'_k = 1 + ha_{k-1}$ satisfies $|hA'| = |hA|$.
\et

\begin{proof}
If $a_k > 1 + ha_{k-1}$, then there is a positive integer $\delta_{k-1}$ such that 
\begin{align*} 
a_k - a_{k-1} 
& = 1 + \delta_{k-1} +  (h-1)a_{k-1} \\ 
& = 1 + \delta_{k-1} +  (h-1)\max(a_{k-1}-1,a_k-a_k). 
\end{align*}  
The construction in Theorem~\ref{compression:theorem:Delta-1} 
produces the set $A' = \{a_1 < \cdots < a_{k-1} < a'_k\}$ with 
\[
a'_k = a_k - \delta_{k-1} = 1 + ha_k 
\] 
and $|hA'| = |hA|$.  This  completes the proof.  
\end{proof}

\bt                                           \label{compression:theorem:Delta}
Let $k \geq 3$ and let $A = \{a_1 <  a_2 < \cdots < a_k\}$ be a set of $k$ integers. 
For every positive integer $h$, there is a set $A' = \{a'_1 <  a'_2 < \cdots < a'_k\}$ 
such that $|hA'| = |hA|$ and 
\beq                                     \label{compression:Delta-ineq-x}
a'_{j+1} - a'_j \leq  1+ (h-1)\max (a'_j - a'_1, a'_k - a'_{j+1})  
\eeq
for all $j \in [1,k-1]$. 
\et

\begin{proof}
If the set $A$ does not satisfy condition~\eqref{compression:Delta-ineq-2} 
for some $j \in [1,k-1]$, then apply Theorem~\ref{compression:theorem:Delta-1} 
and obtain a set $A'$ with $|hA'| = |hA|$ and $\diam(A') < \diam(A)$.  
Repeat  this argument if the set $A'$ does not satisfy inequality~\eqref{compression:Delta-ineq-x} 
for some $j \in [1, k-1]$.  
Continuing this procedure, we obtain a sequence 
of sets whose diameters
form a strictly decreasing sequence of positive integers, 
and so the process must stop after finitely many iterations. 
This  completes the proof.  
\end{proof} 

Theorem~\ref{compression:theorem:Delta-ak} is a special case 
of the following result. 

\bt                                        \label{compression:theorem:ShortForm}
Let $k \geq 3$ and let $A = \{a_1 <  a_2 < \cdots < a_k\}$ be a set of $k$ integers 
such that $a_1=0$ and  
\beq              \label{compression:ShortForm}
a_{j+1} - a_j \leq 1 +  (h-1)\max (a_j  , a_k - a_{j+1})  
\eeq
for all $j \in [1,h-1]$, 
There is a unique integer  $r$  in $[2,k-1]$ such that 
\[
a_{r-1} + a_r <  a_k \leq a_r + a_{r+1}.
\]
For all $j \in [r,k-1 ]$,
\[
a_k \leq 1 + h + h^2 + \cdots + h^{k-j-1} + h^{k-j} a_j. 
\]
\et

\begin{proof}
Because $a_1=0$ and $k \geq 3$, we have 
\[
 a_1+a_2 = a_2  <  a_k < a_{k-1} + a_k.
\]
Choose $r \in [2,k-1]$ such that 
\[
 a_{r-1} + a_r < a_k \leq a_r + a_{r+1}.
\]
For all $j \in [r,k-1 ]$, we have 
\[
  a_k \leq a_r + a_{r+1}  \leq a_j + a_{j+1}
\]
and so 
\[
a_k - a_{j+1} \leq a_j.
\]
It follows that 
\[
a_{j+1} - a_j \leq 1 +  (h-1)\max (a_j,  a_k - a_{j+1} ) = 1 + (h-1)  a_j   
\]
and  
\[
a_{j+1} \leq 1 + ha_j.
\]
If $r \leq j \leq j+1 \leq k-1$, then 
\begin{align*}
a_{j+2} & \leq 1 + ha_{j+1} \\
&  \leq 1 + h(1 + ha_j) \\ 
& =  1 + h  + h^2a_j.  
\end{align*}
By induction, if $r  \leq j \leq j + \ell \leq k$, then 
\[
a_{j+\ell}\leq 1 + h + h^2 + \cdots + h^{\ell-1} + h^{\ell} a_j.
\]
Choosing $\ell = k-j$, we obtain 
\[
a_k \leq 1 + h + h^2 + \cdots + h^{k-j-1} + h^{k-j} a_j. 
\]
This completes the proof.
\end{proof}

For example, let $h=3$ and $k=7$.  Consider the set 
\[
A = \{0,2,7,11,70,85,91\}. 
\]
We have $k = |A| = 7$ and $|hA| = |3A| = 80$.  Applying Theorem~\ref{compression:theorem:Delta}, we obtain  
\[
59 = 70 - 11 = 1 + \delta_4 + 2\max(11,91-70) = \delta_4 + 43 
\]
and so $\delta_4 = 16$ and 
\[
A' = \{0,2,7,11,54,69,75\}.
\]
Note that  
\[
43 = 54 - 11 = 1 +   2\max(11,75-54). 
\]
The sumsets $3A$ and $3A'$ are different, that is, $3A \neq 3A'$, but  
\[
|3A| = |3A'| = 80. 
\]
The compressed set $A'$ satisfies condition~\eqref{compression:ShortForm} 
for all $j \in [1,6]$.  
We have 
\[
11+54 = 65 < 75 < 123 = 54 + 69.
\]
Applying Theorem~\ref{compression:theorem:ShortForm} with $r = 5$, we obtain 
\begin{align*}
75 & < 1 + 3 + 3^2 \cdot 54 \\
75 & < 1 + 3  \cdot 69.
 \end{align*}

\section{Additive isomorphism $\mcr_{\Z^n}(h,k) = \mcr_{\Z}(h,k)$}

Let $X$ and $Y$ be additive abelian groups.  Let $A$ be a finite subset 
of $X$ and $B$ a finite subset of $Y$. 
The function $f:A\rightarrow B$ is a \emph{Freiman isomorphism of order $h$} 
if $f$ is a bijection and, for all $a_1,\ldots, a_h ,a'_1,\ldots, a'_h \in A$, we have 
\[
a_1 + \cdots + a_h = a'_1 + \cdots + a'_h
\]
if and only if  
\[
f( a_1) + \cdots +  f( a_h) = f( a'_1) + \cdots + f( a'_h).
\] 
Thus, $f$ induces a bijection from $hA$ to $hB$ and so $|hA| = |hB|$. 
For applications of Freiman isomorphisms, see~\cite{nath96bb}.

The $\ell^{\infty}$-norm of a vector $\mbx = (x_1,\ldots,x_n) \in \R^n$ is 
\[
\|\mbx\|_{\infty} = \max( |x_j|:j \in [1,n]\}.
\]

The proof of the following result is a simple application of 
Freiman isomorphisms.  

\bt        \label{compression:theorem:LatticePoints}
For all positive integers $h$, $k$, and $n$, 
\[
\mcr_{\Z^n}(h,k) = \mcr_{\Z}(h,k). 
\]
\et

\begin{proof}
The group $\Z^n$ contains  subgroups  isomorphic to \Z, 
and so $\mcr_{\Z}(h,k) \subseteq  \mcr_{\Z^n}(h,k)$.  
To prove that $\mcr_{\Z^n}(h,k) \subseteq  \mcr_{\Z}(h,k)$, 
it suffices to show that, for every positive integer $h$ and every finite 
set $A$ of lattice points, there is a finite set $B$ of integers that is 
Freiman isomorphic of order $h$ to $A$.

Let $A \subseteq \Z^n$ with $|A|=k$.  
By translation, we can assume that $A \subseteq \N_0^n$.  
Choose a positive integer 
\[
g > h \max\{\|\mba  \|_{\infty}: \mba  \in A\}. 
\]
Define $f:\Z^n \rightarrow \Z$ by 
\[
f(a_1,\ldots, a_n) = \sum_{j=1}^n a_j g^{j-1} 
\]
and let 
\[
B = f(A) \subseteq \N_0.
\]  
Let $\mba_i = (a_{i,1},\ldots, a_{i,n}) \in A$ and  
$\mba'_i = (a'_{i,1},\ldots, a'_{i,n}) \in A$ 
for all $i \in [1, h]$. 
We have 
\[
\sum_{i=1}^{h} \mba_i = \sum_{i=1}^{h}  (a_{i,1},\ldots, a_{i,n}) = \left( \sum_{i=1}^{h}  a_{i,1},\ldots, \sum_{i=1}^{h}  a_{i,n} \right) 
\]
and 
\[
 \sum_{i=1}^{h} \mba'_i  = \sum_{i=1}^{h} (a'_{i,1},\ldots, a'_{i,n})
=   \left(  \sum_{i=1}^{h} a'_{i,1},\ldots,  \sum_{i=1}^{h} a'_{i,n} \right). 
\]
Moreover, 
\[
f(\mba_i) = \sum_{j=1}^n a_{i,j}g^{j-1} 
\qqand 
f(\mba'_i) = \sum_{j=1}^n a'_{i,j}g^{j-1}.
\] 
We have 
\[
 \sum_{i=1}^{h} f(\mba_i) =  \sum_{i=1}^{h} \sum_{j=1}^n a_{i,j} g^{j-1} 
 = \sum_{j=1}^n \sum_{i=1}^h a_{i,j}g^{j-1} 
\] 
and, for all $j \in [1,n]$,  
\[
0 \leq \sum_{i=1}^h a_{i,j} \leq \sum_{i=1}^h  \|\mba_i \|_{\infty} 
\leq  h \max\{\|\mba  \|_{\infty}: \mba  \in A\} < g.
\]
Similarly, 
\[
 \sum_{i=1}^{h} f(\mba'_i) =  \sum_{i=1}^{h} \sum_{j=1}^n a'_{i,j} g^{j-1} 
 = \sum_{j=1}^n \sum_{i=1}^h a'_{i,j}g^{j-1} 
\] 
and, for all $j \in [1,n]$,   
\[
0 \leq \sum_{i=1}^h a'_{i,j} \leq \sum_{i=1}^h  \|\mba'_i \|_{\infty} 
\leq  h \max\{\|\mba'  \|_{\infty}: \mba'  \in A\} < g.
\]
It follows from the uniqueness of the $g$-adic representation that 
\[
 \sum_{i=1}^{h} f(\mba_i) =  \sum_{i=1}^{h} f(\mba'_i) 
\]
if and only if, for all $j \in [1,n]$,   
\[
\sum_{i=1}^h a_{i,j} =  \sum_{i=1}^h a'_{i,j} 
 \]
if and only if 
\[
\sum_{i=1}^{h} \mba_i = \sum_{i=1}^{h} \mba'_i 
\]
and so the sets $A \subseteq \Z^n$ and $B \subseteq \Z$ are Freiman isomorphic of order $h$.  
This completes the proof. 
\end{proof}

\section{Compression theorem for sums of sets of lattice points}

Let $(\mbx,\mby)$ be an inner product on $\R^n$ and let $\mbu$ be a unit vector in $\R^n$.  
For all $r   \in \R$, we have the hyperplane
\[
H_{\mbu}(r) = \{\mbx \in \R^n: (\mbu,\mbx) = r \}.  
\]
The \emph{$(s,t)$-slab} of \mbu\ in $\R^n$  is the convex set 
\beq           \label{compression:slab}
H_{\mbu}(s,t) = \{\mbx \in \R^n : s \leq (\mbu,\mbx) \leq t \}. 
\eeq
We have $H_{\mbu}(s,t) = \emptyset$ if $s > t$.  
If $r < s$, then 
\beq                                 \label{compression:H-disjoint}
H_{\mbu}(q,r) \cap H_{\mbu}(s,t) = \emptyset.
\eeq
For $q\leq r$ and $s \leq t$, there is the additive relation  
\beq                                 \label{compression:H-add}
H_{\mbu}(q,r)+H_{\mbu}(s,t) \subseteq H_{\mbu}(q+s,r+t).
\eeq

The \emph{diameter} with respect to the unit vector \mbu\ 
of a nonempty subset $A$ of $\R^n$ is
\[
\diam_{\mbu}(A) = \sup\{(\mba,\mbu)-(\mba',\mbu):\mba,\mba' \in A\}.
\]

\bt                              \label{compression:theorem:big}
Let $A$ be a finite subset of $\R^n$ with $|A| = k$.  Let $\mbu$ be a unit vector in  $\R^n$.  
Enumerate  the vectors in $A = \{\mba_1,  \mba_2, \ldots, \mba_k\}$ so that 
\beq             \label{compression:increasing}
(\mbu,\mba_1) \leq (\mbu,\mba_2 ) \leq \cdots \leq (\mbu,\mba_k).
\eeq
Let  $h \geq 2$ and suppose that, for some $j \in [1,k-1]$, 
there exist $\alpha_j > 0$ and $\delta_j > 0$ such that 
\beq              \label{compression:lattice-ineq-1}
(\mbu,  \mba_{j+1} - \mba_j) = \alpha_j + \delta_j 
+ (h-1)\max\left((\mbu,\mba_j - \mba_1), (\mbu,\mba_k- \mba_{j+1}) \right).
\eeq
The integer $j$ is unique.  
There is a subset $A' = \{\mba'_1,\ldots, \mba'_k\}$ of $\R^n$ 
with $|A'| = k$ and 
\[
(\mbu,\mba'_1) \leq (\mbu,\mba'_2 ) \leq \cdots \leq (\mbu,\mba'_k)
\]
such that $|hA| =  |hA'|$ and 
\beq             \label{compression:big}
(\mbu,\mba'_{j+1} - \mba'_j) 
= \alpha_j + (h-1) \max( (\mbu, \mba'_j-\mba'_1), (\mbu,\mba'_k - \mba'_{j+1}))
\eeq
 Moreover, 
 \[
\diam_{\mbu}(A') = \diam_{\mbu}(A) - \delta_j. 
 \]
\et

\begin{proof}
Suppose there exist distinct integers $j_1 \in [1,k-1]$ and $j_2\in [1,k-1]$ 
that satisfy~\eqref{compression:lattice-ineq-1}.  
If $j_1 < j_2$, then 
\begin{align*}  
(\mbu, \mba_{j_1+1} - \mba_{j_1} )
& >  (h-1)\max ((\mbu, \mba_{j_1} - \mba_1), (\mbu, \mba_k - \mba_{j_1+1}) ) \\
& \geq  (h-1) (\mbu, \mba_k - \mba_{j_1+1}) \\
&  \geq  (\mbu, \mba_{j_2+1} - \mba_{j_2} ) \\
& >   (h-1)\max ((\mbu, \mba_{j_2} - \mba_1) , (\mbu, \mba_k - \mba_{j_2+1}) ) \\
& \geq   (h-1) (\mbu,\mba_{j_2} - \mba_1)  \\
&  \geq (\mbu, \mba_{j_1+1} - \mba_{j_1}  )
\end{align*}
which is absurd.  Therefore, at most one integer $j \in [1,k-1]$ satisfies~\eqref{compression:lattice-ineq-1}.

Suppose that relation~\eqref{compression:lattice-ineq-1} holds 
 for some $j \in [1,h-1]$.   Let 
\begin{align*}
q =  (\mbu,\mba_1) \qquad  & \qquad  r =  (\mbu,\mba_j) \\ 
s =  (\mbu,\mba_{j+1})  \qquad & \qquad t =  (\mbu,\mba_k) 
\end{align*}
Then 
\[
q \leq r < s \leq t    
\]
and 
\beq             \label{compression:s-r-ineq}
s-r = \alpha_j + \delta_j + (h-1)\max(r-q, t-s).  
\eeq

Let  
\[
A = B \cup C     
\]
where  
\[
B = \{\mba_1,\ldots, \mba_j\} \subseteq H_{\mbu}(q,r)
\] 
and 
\[
C = \{\mba_{j+1},\ldots, \mba_k \} \subseteq H_{\mbu}(s,t) 
\]
and  $H_{\mbu}(q,r)$ and $H_{\mbu}(s,t)$ are the slabs defined by~\eqref{compression:slab}.
We have  
\[
hA = \bigcup_{i=0}^h \left( (h-i)B + iC \right). 
\]
The additive relation~\eqref{compression:H-add} implies that   
\[
(h-i)B \subseteq H_{\mbu}((h-i)q,(h-i)r) 
\]
and 
\[
iC \subseteq H_{\mbu}(is,it) 
\]
and so 
\[
(h-i)B + iC \subseteq H_{\mbu}( (h-i)q + is, (h-i)r + it). 
\]
Let $i \in [0,h-1]$.  If 
\[
 (h-i)r + it <  (h-i-1)q + (i+1)s
\]
or, equivalently, if 
\[
 (h-1-i)(r-q) + i(t-s) < s - r
\]
then, by~\eqref{compression:H-disjoint}, the sets $(h-i)B + iC$ and $(h-1-i)B + (i+1)C$ are disjoint.

From~\eqref{compression:s-r-ineq} we have 
\beq                       \label{compression:qrst-ineq}
 (h-1-i)(r-q) + i(t-s) \leq (h-1)\max(r-q, t-s) < s-r   
\eeq 
and so the sets $(h-i)B + iC$  are pairwise disjoint 
for all $i \in [0,h]$.  It follows that   
\[
|hA| = \sum_{i=0}^h \left| (h-i)B + iC \right|. 
\]

Let $\mba'_{\ell} = \mba_{\ell} $ for all $\ell \in [1,j]$ 
and let 
\[
B' =  \{\mba'_{1},\ldots, \mba'_j\} = B.
\]Then 
\[
B' \subseteq H_{\mbu}(q,r).
\]
Let $\mba'_{\ell} = \mba_{\ell} -\delta_j \mbu$ for all $\ell \in [j+1,k]$ 
and let 
\[
C' =  \{\mba'_{j+1},\ldots, \mba'_k\} = C - \delta_j \mbu .
\]
Because $\mbu$ is a unit vector, we have $\delta_j = (\mbu,\delta_j \mbu)$. 
For all $\ell \in [j+1,k]$, we have 
\begin{align*}
s-\delta_j  & = (\mbu,\mba_{j+1})  -  \delta_j  = (\mbu,\mba_{j+1})  -  (\mbu,\ \delta_j \mbu) \\
& = (\mbu,\mba_{j+1} -\delta_j \mbu) = (\mbu,\mba'_{j+1}) \\
&   \leq (\mbu,\mba_{\ell}) - \delta_j = (\mbu,\mba'_{\ell})  \\
& \leq  (\mbu,\mba_k) - \delta_j = (\mbu,\mba'_k) \\ 
& = t-\delta_j
\end{align*}
and so 
\[
C' \subseteq H_{\mbu}(s-\delta_j,t-\delta_j).
\]

Let $A' = B' \cup C'$.  We have 
\[
hA' = \bigcup_{i=0}^h \left( (h-i)B' + iC' \right) 
\]
and 
\[
 (h-i)B' + iC' \subseteq H_{\mbu}((h-i)q+i(s-\delta_j), (h-i)r+i(t-\delta_j))
\]
for all $i \in [0,h]$.  The $h+1$ sumsets $ (h-i)B' + iC'$ are pairwise disjoint if
\[
 (h-i)r+i(t-\delta_j) < (h-1-i)q+(i+1)(s-\delta_j)
\]
or, equivalently, if 
\[
(h-1-i)(r-q)+i(t-s) < s-r-\delta_j
\]
for all $i \in [0,h-1]$.
From relation~\eqref{compression:lattice-ineq-1}, we have 
\begin{align*}
(h-1-i) & (r-q)+i(t-s)  \leq (h-1)\max(r-q,t-s) \\
& =  (h-1)\max\left((\mbu,\mba_j - \mba_1), (\mbu,\mba_k- \mba_{j+1}) \right) \\
& = (\mbu,  \mba_{j+1} - \mba_j) - \alpha_j - \delta_j \\
& < s - r -\delta_j. 
\end{align*} 
Thus,  the sets $(h-i)B' + iC'$ pairwise disjoint for $i \in [0,h]$ and 
\[
|hA'| = \sum_{i=1}^h |(h-i)B' + iC'| .
\]
Because $B'=B$ and $C' = C -\delta_j\mbu$, we have 
\[
(h-i)B' + iC' = (h-i)B  + iC - i\delta_j\mbu
\]
and so 
\[
|(h-i)B' + iC'| = | (h-i)B  + iC - i\delta_j\mbu| = | (h-i)B  + iC| 
\]
and 
\[
|hA| = \sum_{i=0}^h |(h-i)B  + iC| = \sum_{i=0}^h |(h-i)B'  + iC'| = |hA'|.
\]
Moreover,
\[
\diam_{\mbu}(A') = (\mbu,\mba'_k) - (\mbu,\mba'_1) 
 = (\mbu,\mba _k) - (\mbu,\mba _1) - \delta_j = \diam_{\mbu}(A) - \delta_j. 
\]
This completes the proof. 
\end{proof}

\bt                              \label{compression:theorem:bigger}
Let $A = \{\mba_1,\ldots, \mba_k\}$ be a finite subset of $\Z^n$ with $|A| = k$.  
Let $\mbe_1,\ldots, \mbe_n$ be the standard basis for  $\R^n$.  
For all $p \in [1,n]$, there is a permutation $\sigma_p$ of the integer interval  $[1,k]$ such that 
\beq          \label{compression:monotone}
(\mbe_p,\mba_{\sigma_p(1)}) \leq (\mbe_p,\mba_{\sigma_p(2)}) \leq \cdots \leq (\mbe_p,\mba_{\sigma_p(k)}).
\eeq
There exists a set $A^* = \{\mba^*_1,\ldots, \mba^*_k\}$ in $\Z^n$ with $\left|A^*\right|=k$ such that 
$|hA| =  |hA^*|$ and, for all $p \in [1,n]$ and $j \in [1,k-1]$, 
\begin{align*} 
& \left( \mbe_p,\mba^*_{\sigma_p(j+1)} - \mba^*_{\sigma_p(j)} \right) \\
& \qquad \leq 1  +  (h-1) \max\left( \left(\mbe_p, \mba^*_{\sigma_p(j)} - \mba^*_{\sigma_p(1)} \right),
  \left(\mbe_p,\mba^*_{\sigma_p(k)} - \mba^*_{\sigma_p(j+1)} \right) \right). 
\end{align*} 
\et

\begin{proof}
Because the standard basis vectors $\mbe_1,\ldots, \mbe_n$ are pairwise orthogonal, 
the monotonicity inequalities~\eqref{compression:monotone} persist if each vector $\mba_j$ is replaced by 
a vector of the form $\mba_j - \mbv_j$ for any vector $\mbv_j$ in the hyperplane $H_{\mbe_p}(0)$, 
that is, by any vector that is a linear combination of the $n-1$ vectors in the set 
$\{\mbe_1,\ldots, \mbe_n\}\setminus \{\mbe_p\}$. 

If 
\begin{align*} 
& \left( \mbe_1,\mba_{\sigma_1(j+1)} - \mba_{\sigma_1(j)} \right) \\
& \qquad > 1  +  (h-1) \max\left( \left(\mbe_1, \mba_{\sigma_1(j)} - \mba_{\sigma_1(1)} \right),
  \left(\mbe_1,\mba_{\sigma_1(k)} - \mba_{\sigma_1(j+1)} \right) \right)
\end{align*} 
for some $j \in [1,k-1]$, then there is a positive integer $\delta_j$ such that 
\begin{align*} 
& \left( \mbe_1,\mba_{\sigma_1(j+1)} - \mba_{\sigma_1(j)} \right) \\
& \qquad = 1  +  \delta_j + (h-1) \max\left( \left(\mbe_1, \mba_{\sigma_1(j)} - \mba_{\sigma_1(1)} \right),
  \left(\mbe_1,\mba_{\sigma_1(k)} - \mba_{\sigma_1(j+1)} \right) \right). 
\end{align*} 
Applying Theorem~\ref{compression:theorem:big} with the unit vector $\mbe_1$ and $\alpha_j =1$, 
we obtain a set  $A' = \{\mba'_1,\ldots, \mba'_k\}$  in $\Z^n$ with $|A'|=k$ such that $|hA|=|hA'|$ and 
\[
\diam_{\mbe_1}(A') =   \diam_{\mbe_1}(A) - \delta_j \leq \diam_{\mbe_1}(A) - 1.
\]
Moreover, for all $j \in [1,k]$, there is a scalar $c_{1,j}$ such that $\mba'_j = \mba_j - c_{1,j} \mbe_1$. 
If 
\begin{align*} 
& \left( \mbe_1,\mba'_{\sigma_1(j+1)} - \mba'_{\sigma_1(j)} \right) \\
& \qquad > 1  +  (h-1) \max\left( \left(\mbe_1, \mba'_{\sigma_1(j)} - \mba'_{\sigma_1(1)} \right),
  \left(\mbe_1,\mba'_{\sigma_1(k)} - \mba'_{\sigma_1(j+1)} \right) \right)
\end{align*} 
for some $j' \in [1,k-1]$, then we again apply Theorem~\ref{compression:theorem:big}  
obtain a set $A''$  in $\Z^n$ with $|A''|=k$ such that $|hA|=|hA''|$ and 
\[
\diam_{\mbe_1}(A'') \leq \diam_{\mbe_1}(A') - 1 \leq \diam_{\mbe_1}(A) - 2.
\]
Iterating this process, we obtain a sequence of sets whose diameters form a strictly 
decreasing sequence of positive integers.  This process must terminate  
after finitely many iterations and we obtain a set $A^{(1)} = \left\{ \mba_1^{(1)},\ldots, \mba^{(1)}_k\right\}$ 
with $\left| A^{(1)} \right| = k$ such that $|hA^{(1)}| = |hA|$ and 
\begin{align*} 
& \left( \mbe_1,\mba'_{\sigma_1(j+1)} - \mba'_{\sigma_1(j)} \right) \\
& \qquad \leq1  +  (h-1) \max\left( \left(\mbe_1, \mba'_{\sigma_1(j)} - \mba'_{\sigma_1(1)} \right),
  \left(\mbe_1,\mba'_{\sigma_1(k)} - \mba'_{\sigma_1(j+1)} \right) \right) 
\end{align*} 
for all $j \in [1,k-1]$.   Moreover, each vector $\mba_j^{(1)}$  in $A^{(1)}$ differs from 
the vector  $\mba_j$ in $A$ by a scalar multiple of $\mbe_1$. 

Suppose for some $p \in [1,n-1]$ we have constructed a set $A^{(p)}$ in $\Z^n$ with 
$\left| A^{(p)} \right| = k$ such that $|hA^{(p)}| = |hA|$ and, for all $p' \in [1,p]$ 
and for all $j \in [1,k-1]$,  
\begin{align*}  
& \left(\mbe_{p'},  \mba^{\left(p\right)}_{\sigma_{p'}(j+1)} - \mba^{\left(p\right)}_{\sigma_{p'}(j)}   \right) \\
& \qquad \leq 1  + \left(h-1\right) \max\left( \left(\mbe_{p'} ,  
\mba^{\left(p\right)}_{\sigma_{p'}(j)} - \mba^{\left(p\right)}_{\sigma_{p'}(1)} \right),
 \left(\mbe_{p'}, \mba^{\left(p\right)}_{\sigma_{p'}(k)} - \mba^{\left(p\right)}_{\sigma_{p'}(j+1)} \right)\right). 
\end{align*} 
Moreover, each vector in $A^{(p)}$ differs from a vector in $A$ 
by a linear combination of $\mbe_1, \ldots, \mbe_p$. 

We apply Theorem~\ref{compression:theorem:big} to the set $A^{(p)}$ 
with the unit vector $\mbe_{p+1}$ and obtain a set $A^{(p+1)}$ 
with 
$\left| A^{(p+1)} \right| = k$ such that $\left| hA^{(p+1)} \right| = \left|  hA  \right|$ and
\begin{align*}
& \left(\mbe_{p+1},\mba^{\left(p+1\right)}_{j+1} - \mba^{\left(p+1\right)}_j\right) \\
& \quad  \leq 1  + \left(h-1\right) \max\left( \left(\mbe_{p+1}, \mba^{\left(p+1\right)}_j-\mba^{\left(p+1\right)}_1\right),
 \left(\mbe_{p+1},\mba^{\left(p+1 \right)}_k - \mba^{\left(1\right)}_{j+1}\right)\right)
\end{align*}
for all $j \in [1,k-1]$.  
Because each vector in $A^{(p+1)}$ differs from a vector in $A$ 
by a linear combination of $\mbe_1, \ldots, \mbe_{p+1}$, we also have, 
for all $p' \in [1,p]$ and $j \in [1,k-1]$, 
\begin{align*}  
& \left(\mbe_{p'},  \mba^{(p+1)}_{\sigma_{p'}(j+1)} - \mba^{(p+1)}_{\sigma_{p'}(j)}   \right) \\
& \qquad \leq 1  + \left(h-1\right) \max\left( \left(\mbe_{p'} ,  
\mba^{(p+1)}_{\sigma_{p'}(j)} - \mba^{(p+1)}_{\sigma_{p'}(1)} \right),
 \left(\mbe_{p'}, \mba^{(p+1)}_{\sigma_{p'}(k)} - \mba^{(p+1)}_{\sigma_{p'}(j+1)} \right)\right). 
\end{align*} 
Continuing inductively completes the proof. 
\end{proof}

\section{An upper bound for  $N(h,k)$}
The computational complexity of the sumset size problem is determined by the integer 
$N(h,k)$, which is the smallest integer $N$ such that 
\[  
\mcr_{\Z}(h,k) = \{ |hA|: A \subseteq [0,N-1] \text{ and } |A| = k \}.
\]
In this section we compute an upper bound for $N(h,k)$.  
This result is based on work with Noah Kravitz. 

Let $p$ be a prime number and let $(\Z/p\Z)^*$ 
be the multiplicative group of nonzero congruence classes mod $p$. 
For all $a \in \Z$, let $\overline{a}$ denote the congruence class $a+p\Z$. 
For  $\overline{A}\subseteq \Z/p\Z$ and $\lambda \in (\Z/p\Z)^*$, 
let 
\[
\lambda \ast \overline{A} = \{\lambda \overline{a}:\overline{a} \in \overline{A} \}.
\]

The following result is known.  A proof is included for completeness.  	

\bl                    \label{N:lemma}
(a) For  $\overline{A}\subseteq \Z/p\Z$ and $\lambda \in (\Z/p\Z)^*$, 
the sets $\overline{A}$ and $\lambda \ast \overline{A}$ are Freiman isomorphic of order $h$ 
 for  all  $h \geq 2$. 

(b)  Let $A = \{a_1,\ldots, a_k\}$ be a set of $k$ integers and let 
$\overline{A} = \{ \overline{a_1},\ldots, \overline{a_k}  \} $ be the corresponding set of 
congruence classes in the field $\Z/p\Z$. 
Let $h \geq 2$.   If $|a_i| < p/2h$ for all $i \in [1,k]$, 
then the sets $A$ and $\overline{A}$ are Freiman isomorphic of order $h$.
\el

\begin{proof}
To prove (a), it suffices to observe that, 
for all $\lambda \in (\Z/p\Z)^*$, the function $f:\overline{A} \rightarrow \lambda \ast \overline{A}$ 
defined by $f(\overline{a}) = \lambda \overline{a}$ is a Freiman isomorphism of order $h$  
for  all  $h \geq 2$.  

To prove (b), consider the function $g:A \rightarrow \overline{A}$ 
defined by $g(a) = \overline{a}$.  
Because $|a_i| < p/2h < p/2$ for all $i \in [1,k]$, 
 the function $g$ is a bijection.
Let $a_i,b_i \in A$ for all $i \in [1,h]$.  If
\[
\sum_{i=1}^h a_i  = \sum_{i=1}^h b_i  
\]
then 
\[
\sum_{i=1}^h g(a_i)  = \sum_{i=1}^h \overline{a_i} = \sum_{i=1}^h a_i + p\Z  = 
\sum_{i=1}^h b_i + p\Z =  \sum_{i=1}^h \overline{b_i} = \sum_{i=1}^h g(b_i). 
\]
Conversely, if 
\[
\sum_{i=1}^h g(a_i)  = \sum_{i=1}^h g(b_i)
\] 
then 
\[
\sum_{i=1}^h a_i + p\Z  = \sum_{i=1}^h \overline{a_i} = 
 \sum_{i=1}^h \overline{b_i} = \sum_{i=1}^h b_i + p\Z 
\]
and so $p$ divides $\sum_{i=1}^h b_i -  \sum_{i=1}^h a_i$.  
Because $|a_i| < p/2h$ for all $i \in [1,k]$, we have 
 \[
\left|  \sum_{i=1}^h a_i  \right| \leq  \sum_{i=1}^h |a_i | < \frac{p}{2} 
\qqand 
\left|  \sum_{i=1}^h b_i  \right| \leq  \sum_{i=1}^h |b_i |  <\frac{p}{2} 
 \]
and so 
 \[
 -p < \sum_{i=1}^h b_i -  \sum_{i=1}^h a_i  < p.
 \]
 It follows that $ \sum_{i=1}^h a_i =  \sum_{i=1}^h b_i $, 
 and so $g$ is an Freiman isomorphism of order $h$.
 This completes the proof.  
\end{proof}

Let $p$ be an odd prime.  
For $\overline{x} =  x+p\Z\in \Z/p\Z$, let $|\overline{x}|$ be the  smallest  
absolute value of an integer 
in the congruence class $\overline{x}$. Thus, $|\overline{x}| \in [0, (p-1)/2]$. 
For $\mbx = (\overline{x_1},\ldots, \overline{x_k}) \in (\Z/p\Z)^k$, let 
\[
\|\mbx\|_{\infty} = \max\{ | \overline{x_j} |:j \in [1,k]\} \in \left[ 0, \frac{p-1}{2} \right].
\]
For all $\mbv_1, \mbv_2 \in (\Z/p\Z)^k$, we have 
\beq                         \label{compression:l-infinity}
\| \mbv_1 \pm \mbv_2 \|_{\infty} \leq \| \mbv_1 \|_{\infty} + \|\mbv_2 \|_{\infty}.
\eeq

Let $\mbo = (\overline{0},\ldots, \overline{0}) \in  (\Z/p\Z)^k$.  
For $\mbx \in (\Z/p\Z)^k$ and $r \in [0,(p-1)/2]$,  define the boxes 
\[
B_r(\mbo) = \left\{ \mbv \in (\Z/p\Z)^k : \|\mbv\|_{\infty} \leq r  \right\} 
\]
and
\[
B_r(\mbx)  = \left\{   \mbx + \mbv:  \mbv \in B_r(\mbo) \right\} = \mbx +  B_r(\mbo).
\]
We have  
\beq          \label{N:Br}
\left| B_r(\mbo) \right| = \left| B_r(\mbx) \right| =  (2r+1)^k.
\eeq

\bt                  \label{compression:theorem:N(h,k)}
For all $h \geq 3$ and $k \geq 3$,
\[
N(h,k) < 4(8h)^{k-1}.
\]
\et

\begin{proof}
Let $t \in \mcr_{\Z}(h,k)$ and let $M = M(t)$ be the smallest positive integer 
such that there exists a set 
\[
A = \{a_1,\ldots, a_k\} \subseteq [-M,M] 
\]
with $|A|=k$ and $|hA|=t$.  
We shall prove that $M(t) < 2(8h)^{k-1}$ for all $t \in \mcr_{\Z}(h,k)$.  
It follows that $N(h,k) < 4(8h)^{k-1}$.

By Bertrand's postulate, there is a prime $p$ such that   
\beq          \label{N:1} 
2h M < p < 4hM   
\eeq
 and so  
\[
|a_j| \leq M < \frac{p}{2h} 
\]
for all $j \in [1,k]$.   
Let 
\[
\overline{A} = \{\overline{a_1},\ldots, \overline{a_k} \} \subseteq \Z/p\Z. 
\]
By Lemma~\ref{N:lemma}(b), 
the sets $A$ and $\overline{A}$ are Freiman isomorphic of order $h$. 
By Lemma~\ref{N:lemma}(a), the sets $\overline{A}$ and $\lambda \ast \overline{A}$ 
are Freiman isomorphic of order $h$ for all $\lambda \in (\Z/p\Z)^*$.  

Choose the integer $r$ such that  
\[
 \frac{1}{2} \left( p^{1-(1/k)} -1  \right) < r \leq  \frac{1}{2} \left( p^{1-(1/k)} +1  \right). 
\]
Then
\beq          \label{N:2}
p(2r+1)^k > p^k.
\eeq
We may assume that $M \geq 2(4h)^{k-1}$ and so  
\[
(4h)^k \leq 2hM < p.
\]
It follows that  
\beq          \label{N:r}
2r \leq  p^{1-(1/k)} + 1 < \frac{ 2p}{p^{1/k}}   <   \frac{p}{2h}.
\eeq

Fix an ordering of the $k$ congruence classes in the set $\overline{A} $  
and consider the vector 
\[
\mba = ( \overline{a_1},\ldots, \overline{a_k} ) \in (\Z/p\Z)^k.
\]
For all $\lambda  \in  \Z/p\Z $,  we have the vector 
\[
\lambda \mba =  (\lambda \overline{a_1},\ldots, \lambda \overline{a_k} ) 
\]
and   the box $B_r(\lambda \mba)$ in $(\Z/p\Z)^k$.  
From relations~\eqref{N:Br} and~\eqref{N:2}, we obtain 
\[
\sum_{\lambda \in  \Z/p\Z } \left| B_r(\lambda \mba) \right| 
= p(2r+1)^k > p^k = \left|(\Z/p\Z)^k\right| 
\]
and so the $p$ boxes $B_r(\lambda \mba)$ are not pairwise disjoint.  
It follows that there exist $\lambda_1, \lambda_2 \in \Z/p\Z$ such that 
$\lambda_1 \neq \lambda_2$ and 
\[
B_r(\lambda_1 \mba) \cap B_r(\lambda_2 \mba) \neq \emptyset. 
\]
We obtain $\mbv_1,\mbv_2 \in B_r(\mbo)$ such that
\[
\lambda_1 \mba + \mbv_1 = \lambda_2 \mba + \mbv_2.
\]
Let 
\[
\lambda = \lambda_1-\lambda_2 \in (\Z/p\Z)^* 
\]
and  
\[
\mbb = \lambda \mba  = \lambda_1 \mba   - \lambda_2 \mba  = \mbv_2 - \mbv_1 
\in (\Z/p\Z)^k. 
\] 
From inequalities~\eqref{compression:l-infinity} and ~\eqref{N:r}, we have  
\[
\|\mbb\|_{\infty} = \| \mbv_2 - \mbv_1\|_{\infty} 
 \leq  \| \mbv_1 \|_{\infty}+ \| \mbv_2\|_{\infty} \leq 2r  < \frac{p}{2h}
\]
and so there are integers $b_1,\ldots, b_k$ such that 
\[
\mbb =  ( \overline{b_1},\ldots, \overline{b_k} ) = \lambda \mba
\]
and 
\[
|b_j| \leq 2r < \frac{p}{2h}
\]
for all $j \in [1,h]$.  
Let 
\[
\overline{B} = \{ \overline{b_1},\ldots, \overline{b_k} \} = \lambda \ast \overline{A}. 
\] 
By Lemma~\ref{N:lemma}, the sets of congruence classes $\overline{A} $ 
and $\overline{B}$ are Freiman isomorphic of order $h$, 
and   the set  of congruence classes $\overline{B} $ 
and the set of integers $B = \{b_1,\ldots, b_k\}$ are Freiman isomorphic of order $h$.   
Therefore,  $|hB| = |h\overline{B}| = |h\overline{A}| = |hA| = t$.  
Moreover, 
\[
B \subseteq [-2r,2r]. 
\]
The minimality of $M$ implies 
\begin{align*}
M & \leq 2r \leq 2 p^{1-(1/k)}  < 2(4hM)^{1-(1/k)}. 
\end{align*}
Solving for $M$, we obtain 
\[
M   < 2(8h)^{k-1}.
\]
This completes the proof.  
\end{proof}

\end{document}